\documentclass[11pt,a4paper,twoside]{amsart}

\usepackage{graphics}
\usepackage{color}
\usepackage{epic}
\usepackage{mfpic}

\newtheorem{theo}{\bf Theorem}[section]

\newtheorem{lemma}[theo]{\bf Lemma}

\newtheorem{coro}[theo]{\bf Corollary}

\newcommand{\Z}{\mathbb Z}

\newcommand{\uarrow}{{\begin{picture}(0,0) 
      \put(0,0){\line(5,6){5}}
      \put(3.33,4){\vector(3,4){0}}
    \end{picture}}}

\newcommand{\darrow}{{\begin{picture}(0,0)
      \put(0,0){\line(5,-4){5}}
      \put(3.33,-2.667){\vector(4,-3){0}}
    \end{picture}}}

\newcommand{\udashline}{{\begin{picture}(0,0)
      \dottedline{0.5}(0,0)(5,6)
    \end{picture}}}

\newcommand{\ddashline}{{\begin{picture}(0,0)
      \dottedline{0.5}(0,0)(5,-4)
    \end{picture}}}

\newcommand{\udigitline}[1]{{\begin{picture}(0,0)
      \dottedline{0.5}(0,0)(1.5,1.8)
      \dottedline{0.5}(3.5,4.2)(5,6)
      \put(1.5,2){#1}
    \end{picture}}}

\newcommand{\ddigitline}[1]{{\begin{picture}(0,0)
      \dottedline{0.5}(0,0)(1.5,-1.2)
      \dottedline{0.5}(3.5,-2.8)(5,-4)
      \put(1.5,-3){#1}
    \end{picture}}}

\begin{document}

\title{Cylindrical lattice paths and the Loehr-Warrington
  $10^n$ conjecture}
\author{Jonas Sj{\"o}strand}
\address{Department of Mathematics, Royal Institute of Technology \\
  SE-100 44 Stockholm, Sweden}
\email{jonass@kth.se}
\keywords{lattice path; cylinder graph; weight function}
\subjclass{Primary: 05A15; Secondary: 05C38}
\date{29 September 2005}

\begin{abstract}
  The following special case of a conjecture by Loehr and Warrington
  was proved recently by Ekhad, Vatter, and Zeilberger:

  There are $10^n$ zero-sum words of length $5n$ in the alphabet $\{+3,-2\}$
  such that no zero-sum consecutive subword that starts with
  $+3$ may be followed immediately by $-2$.

  We give a simple bijective proof of the conjecture in its original
  and more general setting. To do this we reformulate the problem
  in terms of cylindrical lattice paths.
\end{abstract}

\maketitle

\section{Introduction}
\noindent
Let $a$ and $b$ be positive integers.
Given a word $w$ in the alphabet $\{+a,-b\}$,
a zero-sum consecutive subword of $w$ is said to be {\em illegal}
if it starts with $+a$, and $-b$
comes immediately after the subword in $w$.
Example:
$$w=\mbox{}-2+3\underbrace{\mbox{\boldmath $\mbox{}+3$}-2-2-2+3}
_{\mbox{\small illegal subword}}\mbox{\boldmath $\mbox{}-2$}-2+3$$
We will prove the following:
\begin{theo}\label{th:main}
  If $a$ and $b$ are relatively prime,
  there are ${{a+b}\choose{a}}^n$ zero-sum words of length $(a+b)n$
  in the alphabet $\{+a,-b\}$ without illegal subwords.
\end{theo}

Example: If $a=2$, $b=1$ and $n=2$, the ${{2+1}\choose2}^2=9$
words counted in the theorem are
$$\mbox{}+2-1-1+2-1-1$$
$$\mbox{}-1+2-1+2-1-1$$
$$\mbox{}-1-1+2+2-1-1$$
$$\mbox{}-1-1-1+2+2-1$$
$$\mbox{}-1+2-1-1+2-1$$
$$\mbox{}-1-1+2-1+2-1$$
$$\mbox{}-1-1-1+2-1+2$$
$$\mbox{}-1-1-1-1+2+2$$
$$\mbox{}-1-1+2-1-1+2$$

The theorem was conjectured by Nick Loehr and Greg Warrington.
In a recent paper~\cite{zeilberger} Shalosh Ekhad, Vince Vatter, and
Doron Zeilberger proved the special case $a=3$, $b=2$, using
a computer. Inspired by their proof, Loehr and Warrington, together
with Bruce Sagan~\cite{sagan},
found a computer-free proof of the more general
case when $b=2$ and $a$ is any odd positive integer.
We greatly admire the automatical approach of Ekhad, but we feel
that a beautiful problem like this ought to have a beautiful
solution. And indeed it has!

First we will present a geometrical construction where the words in
the alphabet $\{+a,-b\}$ are interpreted as paths on a
cylinder graph. Then we will give bijections between
these paths, certain weight functions
on the edges of the graph, and ordered sequences of
lap cycles. The latter ones are easy to count.
Finally, in the last section we examine what happens if $a$ and $b$ have
a common factor.

\section{The geometrical construction}
In the following we let $a$ and $b$ be any positive integers.

After thinking about the Loehr-Warrington conjecture for a while,
most people will probably discover the following natural
reformulation:

You live in a skyscraper $\Z$. In the morning you get your exercises
by climbing out through the window, following $(a+b)n$ one-way
ladders, and climbing into your apartment again. At each level there is
one ladder going $a$ levels up and another ladder going $b$
levels down.
Once you have climbed up from a level you never
climb down from that level anymore that morning.
In how many ways can you perform your exercises?

Now here is the key observation:
Since $a\equiv-b\pmod{a+b}$, after $x$
ladders we are at a level $y$ such that
$y\equiv ax\pmod{a+b}$.
We define a directed graph $G_{a,b}$ whose vertex set is
the subset of the infinite cylinder $\Z_{a+b}\times\Z$
consisting of all points $(x,y)$ such that $y\equiv ax\pmod{a+b}$.
From every vertex point $(x,y)$ there is an up-edge
$(x,y)\rightarrow(x+1,y+a)$ and a down-edge
$(x,y)\rightarrow(x+1,y-b)$.
If $a$ and $b$ are relatively prime, no two points in $G_{a,b}$ have
the same $y$-coordinate\footnote{This is the only time we use the assumption
in Theorem~\ref{th:main}
that $a$ and $b$ are relatively prime.\label{fn:prime}}.
We have mapped the ladders
to the cylinder such that no ladders intersect!

\begin{figure}
\begin{center}
  \setlength{\unitlength}{1mm}
  \begin{picture}(121,62)(-21,5)
    \newsavebox{\Gexempel}
    \savebox{\Gexempel}(0,0)[lb]{\setlength{\unitlength}{0.9mm}%
      \begin{picture}(50,70)(-12,-50)
    \put(0,0){\circle*{2}}
    \put(25,0){\circle*{2}}
    \put(-11,-1){$(0,0)$}
    \put(27,-1){$(0,0)$}

    \multiput(0,10)(5,6){2}{\circle*{1}}
    \multiput(0,0)(5,6){4}{\circle*{1}}
    \multiput(0,-10)(5,6){5}{\circle*{1}}
    \multiput(0,-20)(5,6){6}{\circle*{1}}
    \multiput(5,-24)(5,6){5}{\circle*{1}}
    \multiput(15,-22)(5,6){3}{\circle*{1}}
    \multiput(20,-26)(5,6){2}{\circle*{1}}

    \multiput(0,10)(5,6){1}{\uarrow}
    \multiput(0,0)(5,6){3}{\uarrow}
    \multiput(0,-10)(5,6){4}{\uarrow}
    \multiput(0,-20)(5,6){5}{\uarrow}
    \multiput(5,-24)(5,6){4}{\uarrow}
    \multiput(15,-22)(5,6){2}{\uarrow}
    \multiput(20,-26)(5,6){1}{\uarrow}

    \multiput(0,10)(5,6){2}{\darrow}
    \multiput(0,0)(5,6){4}{\darrow}
    \multiput(0,-10)(5,6){5}{\darrow}
    \multiput(0,-20)(5,6){5}{\darrow}
    \multiput(10,-18)(5,6){3}{\darrow}
    \multiput(15,-22)(5,6){2}{\darrow}
    
    \dashline[50]{2}(0,-30)(0,20)
    \dashline[50]{2}(25,-30)(25,20)
    \put(5,-33){\Huge$\vdots$}
    \put(11.5,-31){\Huge$\vdots$}
    \put(18,-33){\Huge$\vdots$}

    \put(5,20){\Huge$\vdots$}
    \put(11.5,22){\Huge$\vdots$}
    \put(18,20){\Huge$\vdots$}

  \end{picture}}

\put(0,-5){\usebox{\Gexempel}}
\put(47,35){\Huge$\Leftrightarrow$}
\put(60,0){\resizebox{40mm}{!}{\includegraphics{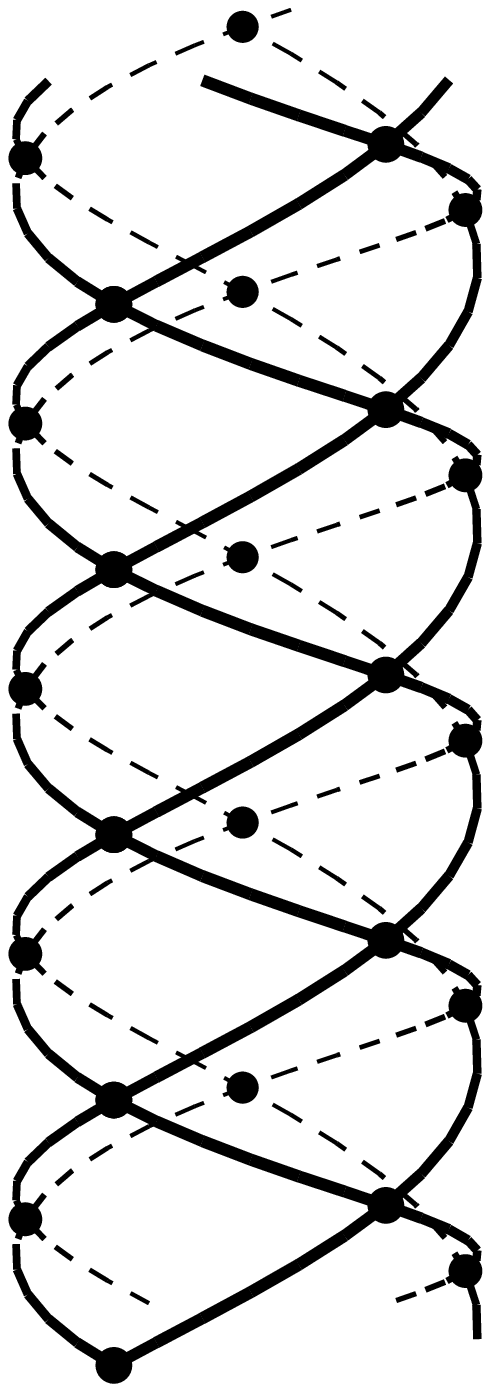}}}
\put(60,38){$(0,0)$}
\put(69,39){\vector(1,0){4}}

\put(-20,30){\vector(1,0){10}}
\put(-20,30){\vector(0,1){10}}
\put(-8,29){$x$}
\put(-21,43){$y$}

\end{picture}
\end{center}
\caption{The graph $G_{3,2}$ represented as a vertical strip to the left.
(We have stretched the $x$-axis by a factor $\sqrt{ab}$ to
make the lattice rectangular, but this is merely cosmetics.)
When the borders are welded together the result is the cylinder to
the right.\label{fig:cylinder}}
\end{figure}

Figure~\ref{fig:cylinder} shows a graphical representation
of $G_{3,2}$. It is an infinite vertical strip whose borders
are welded together. The points with $x=0$
constitute the {\em weld} and are called {\em weld points}.

As far as we know, no one has studied this graph before. The closest
related research we could find is two papers about nonintersecting
lattice paths on the cylinder, one by Peter Forrester~\cite{forrester}
and one by Markus Fulmek~\cite{fulmek}. Curiously, their paths 
essentially go along the axis of the cylinder while ours
essentially go around it!

Let us fix the following graph terminology: A {\em path} is an
ordered sequence of vertices $v_0v_1,\ldots,v_m$ such that
there is an edge from $v_{i-1}$ to $v_i$ for $i=1,2,\ldots,m$.
(Repeated vertices and edges are allowed.) The integer $m$ is the {\em length}
of the path. If $v_0=v_m$ the path is called a {\em cycle}.

A path on $G_{a,b}$ may leave a certain vertex several times, sometimes
going down, sometimes going up. If for each vertex all downs come before
all ups, the path is said to be {\em downs-first}.
Now Theorem~\ref{th:main} can be reformulated:
\begin{equation}\label{eq:G}
  \parbox{100mm}{\em
    There are ${{a+b}\choose a}^n$ downs-first cycles on $G_{a,b}$ of length
    $(a+b)n$ starting at the origin.
  }
\end{equation}
We will prove~(\ref{eq:G}) for any positive integers $a$ and $b$.
(Note that this does not imply that Theorem~\ref{th:main} is
true if $a$ and $b$ have a common factor, see footnote~\ref{fn:prime}.)

A path of length $a+b$ that starts and ends on the weld is called
a {\em lap}, and a lap that is a cycle is called a {\em lap cycle}.
Obviously, there are ${{a+b}\choose a}$ different lap cycles starting
at the origin. Our proof will be a bijection that maps downs-first
cycles to a sequence of lap cycles. The following lemma is crucial.

\begin{lemma}\label{lm:downs-first}
  A downs-first cycle beginning at a weld point never
  visits higher weld points.
\end{lemma}
\begin{proof}
  Suppose the downs-first cycle, starting at a weld point $p$,
  visits a weld point $q$ higher than $p$.
  Let $pq$ be the path along the cycle from its starting
  point $p$ to $q$ (if the cycle visits $q$ several times,
  choose any visit),
  and let $qp$ be the remaining
  path along the cycle from $q$ to the finish point $p$. Obviously,
  $pq$ and $qp$ must intersect somewhere\footnote{If your cylindrical
    intuition fails you, think like this: The path $pq$ must
    make at least one lap starting at $p$ or below
    and ending at a weld point above $p$. Similarly, the path
    $qp$ must make at least one lap starting at a weld point
    above $p$ and ending at $p$ or below. Clearly these laps
    intersect and cross.}.
  Specifically they must intersect at a point where
  $pq$ goes up and $qp$ goes down.
  But this contradicts the assumption that the cycle is downs-first.
\end{proof}

\begin{figure}
\begin{center}
  \setlength{\unitlength}{0.9mm}
  \begin{picture}(50,70)(-12,-50)
    \put(0,0){\circle*{2}}
    \put(25,0){\circle*{2}}
    \put(-11,-1){$(0,0)$}

    \put(0,-10){\begin{picture}(0,0)
        \multiput(0,10)(5,6){3}{\circle*{1}}
        \multiput(0,0)(5,6){4}{\circle*{1}}
        \multiput(0,-10)(5,6){5}{\circle*{1}}
        \multiput(0,-20)(5,6){6}{\circle*{1}}
        \multiput(5,-24)(5,6){5}{\circle*{1}}
        \multiput(10,-28)(5,6){4}{\circle*{1}}
        \multiput(20,-26)(5,6){2}{\circle*{1}}
        
        \multiput(0,10)(5,6){2}{\uarrow}
        \multiput(0,0)(5,6){3}{\uarrow}
        \multiput(0,-10)(5,6){4}{\uarrow}
        \multiput(0,-20)(5,6){5}{\uarrow}
        \multiput(5,-24)(5,6){4}{\uarrow}
        \multiput(10,-28)(5,6){3}{\uarrow}
        \multiput(20,-26)(5,6){1}{\uarrow}
        
        \multiput(0,10)(5,6){3}{\darrow}
        \multiput(0,0)(5,6){4}{\darrow}
        \multiput(0,-10)(5,6){5}{\darrow}
        \multiput(0,-20)(5,6){5}{\darrow}
        \multiput(5,-24)(5,6){4}{\darrow}
        \multiput(15,-22)(5,6){2}{\darrow}
        
        \dashline[50]{2}(0,-35)(0,20)
        \dashline[50]{2}(25,-35)(25,20)
        \put(5,-37){\Huge$\vdots$}
        \put(11.5,-35){\Huge$\vdots$}
        \put(18,-37){\Huge$\vdots$}
      \end{picture}}
  \end{picture}
\end{center}
\caption{The semi-infinite cylinder graph $H_{3,2}$.\label{fig:H}}
\end{figure}
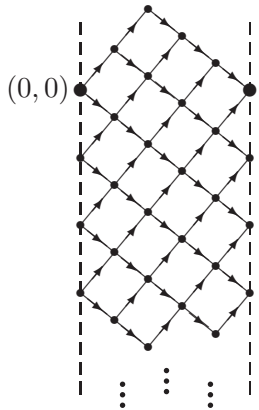

We conclude that the weld points
above the origin, and hence the points above the highest
lap cycle from the origin,
can never be reached by the downs-first cycles
counted in~(\ref{eq:G}).
Let $H_{a,b}$ be the resulting graph
when these points are removed from $G_{a,b}$. Figure~\ref{fig:H}
shows an example.

Now~(\ref{eq:G}) can be slightly reformulated:
\begin{equation}\label{eq:H}
  \parbox{100mm}{\em
    There are ${{a+b}\choose a}^n$ downs-first cycles on $H_{a,b}$ of length
    $(a+b)n$ starting at the origin.
  }
\end{equation}
Before proving~(\ref{eq:H}), and hence our main theorem,
we need some more definitions.

A {\em weight function} on $H_{a,b}$ is an assignment of a
nonnegative integer to every edge in $H_{a,b}$.
The {\em in-weight} and {\em out-weight} of a vertex in $H_{a,b}$
is the sum of the weights of the edges going in to and out from
the vertex, respectively.
A weight function is said to be {\em balanced} if
at each vertex the in-weight and out-weight are equal.
A path in $H_{a,b}$ is said to be {\em covered} by the weight
function if every edge is used by the path at most as many times
as its weight.
The weight function is {\em origin-connected} if for every
vertex with positive out-weight there is a covered path from the origin
to the vertex.

For an example of a balanced origin-connected weight function, see
Figure~\ref{fig:bijections}.

\section{The bijections}
Please keep the cylinder graph $H_{a,b}$ in your mind throughout
the paper.

Since the number of lap cycles beginning at the origin is
${{a+b}\choose a}$, the formulation~(\ref{eq:H}) of our main theorem
follows from the result in this section:
\begin{theo}\label{th:bijections}
  There are bijections between the following three sets:
  \begin{enumerate}
  \item[1.] downs-first cycles of length $(a+b)n$ beginning at the origin,
  \item[2.] balanced origin-connected weight functions with total weight
    sum $(a+b)n$,
  \item[3.] ordered sequences of $n$ lap cycles beginning at the origin.
  \end{enumerate}
\end{theo}
\begin{proof}
  We will define four functions, $f_{1,2}:1\rightarrow2$,
  $f_{2,1}:2\rightarrow1$, $f_{3,2}:3\rightarrow2$, and
  $f_{2,3}:2\rightarrow3$. It should be apparent from the
  presentation below that
  $f_{1,2}\circ f_{2,1}$, $f_{2,1}\circ f_{1,2}$,
  $f_{3,2}\circ f_{2,3}$, and $f_{2,3}\circ f_{3,2}$
  are all identity functions.
  Figure~\ref{fig:bijections} gives an example of the bijections.

\begin{figure}[ht]
\begin{center}
  \setlength{\unitlength}{0.9mm}
  \begin{picture}(145,150)(0,0)

    \newsavebox{\cycle}
    \savebox{\cycle}(0,0)[lb]{\setlength{\unitlength}{0.9mm}%
      \begin{picture}(145,65)(-15,-50)
        \put(-11,-1){$(0,0)$}
        \newsavebox{\temp}
        \savebox{\temp}(0,0)[lb]{%
          \begin{picture}(0,0)
            \put(0,10){\circle*{2}}
            \put(25,10){\circle*{2}}
            \multiput(0,10)(5,6){3}{\circle*{1}}
            \multiput(0,0)(5,6){4}{\circle*{1}}
            \multiput(0,-10)(5,6){5}{\circle*{1}}
            \multiput(0,-20)(5,6){6}{\circle*{1}}
            \multiput(5,-24)(5,6){5}{\circle*{1}}
            \multiput(10,-28)(5,6){4}{\circle*{1}}
            \multiput(20,-26)(5,6){2}{\circle*{1}}
            
            \multiput(0,10)(5,6){2}{\udashline}
            \multiput(0,0)(5,6){3}{\udashline}
            \multiput(0,-10)(5,6){4}{\udashline}
            \multiput(0,-20)(5,6){5}{\udashline}
            \multiput(5,-24)(5,6){4}{\udashline}
            \multiput(10,-28)(5,6){3}{\udashline}
            \multiput(20,-26)(5,6){1}{\udashline}
            
            \multiput(0,10)(5,6){3}{\ddashline}
            \multiput(0,0)(5,6){4}{\ddashline}
            \multiput(0,-10)(5,6){5}{\ddashline}
            \multiput(0,-20)(5,6){5}{\ddashline}
            \multiput(5,-24)(5,6){4}{\ddashline}
            \multiput(15,-22)(5,6){2}{\ddashline}
            
            \dashline[50]{2}(0,-35)(0,20)
            \dashline[50]{2}(25,-35)(25,20)
            \put(5,-37){\Huge$\vdots$}
            \put(11.5,-35){\Huge$\vdots$}
            \put(18,-37){\Huge$\vdots$}
          \end{picture}}

        \multiput(0,-10)(25,0){5}{\usebox{\temp}}
        \sbox{\temp}{}

        \thicklines
        \put(0,0){\darrow}
        \put(5,-4){\uarrow}
        \multiput(10,2)(5,-4){8}{\darrow}
        \put(50,-30){\uarrow}
        \multiput(50,-30)(5,6){4}{\uarrow}
        \put(70,-6){\darrow}
        \multiput(75,-10)(5,6){2}{\uarrow}
        \multiput(85,2)(5,-4){2}{\darrow}
        \put(95,-6){\uarrow}
        \put(100,0){\darrow}
        \multiput(105,-4)(5,6){2}{\uarrow}
        \multiput(115,8)(5,-4){2}{\darrow}
      \end{picture}}

    \put(0,85){\usebox{\cycle}}
    \sbox{\cycle}{}

\newsavebox{\weights}
\savebox{\weights}(0,0)[lb]{\setlength{\unitlength}{1mm}%
  \begin{picture}(50,70)(-12,-50)
    \put(0,0){\circle*{2}}
    \put(25,0){\circle*{2}}
    \put(-11,-1){$(0,0)$}

    \put(0,-10){\begin{picture}(0,0)
        \multiput(0,10)(5,6){3}{\circle*{1}}
        \multiput(0,0)(5,6){4}{\circle*{1}}
        \multiput(0,-10)(5,6){5}{\circle*{1}}
        \multiput(0,-20)(5,6){6}{\circle*{1}}
        \multiput(5,-24)(5,6){5}{\circle*{1}}
        \multiput(10,-28)(5,6){4}{\circle*{1}}
        \multiput(20,-26)(5,6){2}{\circle*{1}}
        
        \dashline[50]{2}(0,-35)(0,20)
        \dashline[50]{2}(25,-35)(25,20)
        \put(5,-37){\Huge$\vdots$}
        \put(11.5,-35){\Huge$\vdots$}
        \put(18,-37){\Huge$\vdots$}
        
        \multiput(0,10)(5,6){2}{\udashline}
        \multiputlist(0,0)(5,6)[lb]{\udigitline1,\udigitline3,\udigitline1}
        \multiput(0,-10)(5,6){4}{\udashline}
        \multiput(0,-20)(5,6){5}{\udigitline1}
        \multiput(5,-24)(5,6){4}{\udashline}
        \multiput(10,-28)(5,6){3}{\udashline}
        \multiput(20,-26)(5,6){1}{\udashline}
        
        \multiputlist(0,10)(5,6)[lb]{\ddigitline2,\ddashline,\ddashline}
        \multiputlist(0,0)(5,6)[lb]{\ddigitline1,\ddashline,%
          \ddigitline2,\ddigitline1}
        \multiputlist(0,-10)(5,6)[lb]{\ddashline,\ddigitline1,\ddashline,%
          \ddigitline2,\ddigitline1}
        \multiputlist(0,-20)(5,6)[lb]{\ddashline,\ddashline,\ddigitline1,%
          \ddashline,\ddigitline2}
        \multiputlist(5,-24)(5,6)[lb]{\ddashline,\ddashline,\ddigitline1,%
          \ddashline}
        \multiputlist(15,-22)(5,6)[lb]{\ddashline,\ddigitline1}
      \end{picture}}
  \end{picture}}

\put(10,0){\usebox{\weights}}
\sbox{\weights}{}

    \newsavebox{\laps}
    \savebox{\laps}(0,0)[lb]{\setlength{\unitlength}{1mm}%
      \begin{picture}(50,70)(-12,-50)
        \put(0,0){\circle*{2}}
        \put(25,0){\circle*{2}}
        \put(-11,-1){$(0,0)$}

        \put(0,-10){\begin{picture}(0,0)
            \multiput(0,10)(5,6){3}{\circle*{1}}
            \multiput(0,0)(5,6){4}{\circle*{1}}
            \multiput(0,-10)(5,6){5}{\circle*{1}}
            \multiput(0,-20)(5,6){6}{\circle*{1}}
            \multiput(5,-24)(5,6){5}{\circle*{1}}
            \multiput(10,-28)(5,6){4}{\circle*{1}}
            \multiput(20,-26)(5,6){2}{\circle*{1}}
            
            \dashline[50]{2}(0,-35)(0,20)
            \dashline[50]{2}(25,-35)(25,20)
            \put(5,-37){\Huge$\vdots$}
            \put(11.5,-35){\Huge$\vdots$}
            \put(18,-37){\Huge$\vdots$}
          \end{picture}}
        
        \put(0.5,0.8){\line(5,-4){4.5}}
        \put(5,-2.8){\line(5,6){9.4}}
        \put(14.4,8.48){\line(5,-4){10}}
        
        \put(0,0){\line(5,-4){5}}
        \put(5,-4){\line(5,6){5.6}}
        \put(10.6,2.72){\line(5,-4){10}}
        \put(20.6,-5.28){\line(5,6){4.4}}
        
        \put(0.6,-10.48){\line(5,6){10}}
        \put(10.6,1.52){\line(5,-4){14.5}}

        \put(0.6,-10.48){\line(5,-4){9.4}}
        \put(10,-18){\line(5,6){9.4}}
        \put(19.4,-6.72){\line(5,-4){5}}

        \put(0.6,-30.48){\line(5,6){9.4}}
        \put(10,-19.2){\line(5,-4){14.5}}

      \end{picture}}

    \put(80,0){\usebox{\laps}}
    \sbox{\laps}{}
    
    \put(35,76){\Huge$\Updownarrow$}
    \put(68,35){\Huge$\Leftrightarrow$}

  \end{picture}
\end{center}
\caption{An example of the bijections in Theorem~\ref{th:bijections}
with $a=3$, $b=2$, and $n=5$. At the top is a downs-first cycle;
below is the corresponding weight function to the left and its packed
sequence of lap cycles to the right.\label{fig:bijections}}
\end{figure}

  \vskip3mm
  $1\rightarrow2$:
  Given a downs-first cycle beginning at the origin,
  to every edge of $H_{a,b}$ we assign a weight that is
  the number of times the edge
  is used by the cycle.
  This weight function is obviously balanced and origin-connected.
  Furthermore, it is the only such function with total sum $(a+b)n$
  that covers the cycle.

  \vskip3mm
  $2\rightarrow1$:
  Given a balanced and origin-connected weight function we construct
  a downs-first cycle as follows: Start at the origin. At each point,
  go down if that edge has positive weight, otherwise go up.
  Decrease the weight of the followed edge by one. Continue until
  you come to a point with zero out-weight. This must be the origin
  so we have created a downs-first cycle $C$.

  We must show that the
  length of $C$ is $(a+b)n$. Suppose not. Then there
  remain some positive weights. Since the original weight function
  was origin-connected there exists a point $p$ on $C$
  with positive out-weight.
  Since the remaining weight function
  is still balanced it covers some cycle $C'$ that contains $p$.
  Now start at $p$ and follow $C$ and $C'$ in parallel until
  $C$ reaches the origin. Since the origin has no remaining
  in- or out-weight, $C'$ must have reached a point on the weld
  below the origin (the other weld points were removed when we constructed
  $H_{a,b}$).
  This implies that $C$ and $C'$
  intersect at a point where $C$ goes up and $C'$ goes down.
  But that is impossible by the construction of $C$.

  Thus $C$ has length $(a+b)n$, and it is easy to see that
  among all downs-first cycles of that length starting at the origin,
  $C$ is the only one that is covered by the given weight function.

  \vskip3mm
  $3\rightarrow2$:
  Given a sequence of lap cycles $C_1,C_2,\ldots,C_n$ starting at
  the origin, translate $C_1,C_2,\ldots,C_{n-1}$ downwards so that,
  for $1\le i\le n-1$, $C_i$ intersects $C_{i+1}$ in at least one point but
  otherwise goes below it.
  Observe that there is a unique way of ``packing'' the cycles
  like that.

  Now let the weight of each edge in $H_{a,b}$ be
  the number of times the edge is used by the cycles.
  The result is obviously a balanced origin-connected weight function.

  \vskip3mm
  $2\rightarrow3$:
  Given a balanced origin-connected weight function,
  by iteration of the following procedure we construct
  $n$ lap cycles. At the beginning of each
  iteration the weight function is always balanced.

  Start at the lowest weld point $p$ with a positive out-weight
  and create a downs-first cycle from there like this:
  In each step, go down if that edge has positive weight, otherwise go up.
  Decrease the weight of the followed edge by one. Stop as soon as
  you reach $p$ again. By Lemma~\ref{lm:downs-first} this downs-first
  cycle never visits another weld point than $p$, which implies
  that it is a lap cycle. (Remember that a path must
  visit the weld every $(a+b)$-th step.)

  After $n$ iterations we have consumed all weights and produced
  a packed sequence of lap cycles $C_1,C_2,\ldots,C_n$.
  Simply translate the lap cycles so that they all start at the origin.

\end{proof}

\section{What if $a$ and $b$ have a common factor?}
The condition that $a$ and $b$ should be relatively prime is not an
essential restriction, as the following corollary to
Theorem~\ref{th:main} shows.
\begin{coro}
  Let $a$ and $b$ be any positive integers, and put $c=\gcd(a,b)$.
  The number of zero-sum words in the alphabet $\{+a,-b\}$
  of length $(a+b)n$ without illegal subwords is
  $${{(a+b)/c}\choose{a/c}}^{cn}.$$
\end{coro}
\begin{proof}
Let $A=a/c$ and $B=b/c$. Clearly the words counted in the corollary are
in one-to-one correspondence with the zero-sum words
in the alphabet $\{+A,-B\}$ of length $(A+B)cn$ without illegal
subwords. According to Theorem~\ref{th:main}
there are ${{A+B}\choose A}^{cn}$ such words.
\end{proof}

However, Theorem~\ref{th:main} can be generalized in a less trivial way:
\begin{theo}\label{th:extra}
  Let $a$ and $b$ be any positive integers.
  There are ${{a+b}\choose a}^n$ zero-sum words in the alphabet $\{+a,-b\}$
  of length $(a+b)n$ without illegal subwords
  whose length is a multiple of $a+b$.
\end{theo}
\begin{proof}
In the proof of Theorem~\ref{th:main} the only time we
make use of the fact that
$a$ and $b$ are relatively prime is when we conclude
that no two points in $G_{a,b}$ have
the same $y$-coordinate. We need this in order to show that
the vertices in $G_{a,b}$ are
in one-to-one correspondence with the levels in the skyscraper.

However, even if $a$ and $b$ have a common factor, the condition
that the word should have no illegal subword
whose length is a multiple of $a+b$
is equivalent to the downs-first condition on the corresponding
path on $G_{a,b}$. (It just happens that the phrase
``whose length is a multiple of $a+b$'' is unnecessary if
$a$ and $b$ are relatively prime.) Thus we can simply bypass the
skyscraper nonsense and go directly to the proposition~(\ref{eq:G})
which is still valid.
\end{proof}

\end{document}